\documentclass{amsart}
\usepackage{amssymb,amsfonts,amscd,amstext}
\usepackage{epsfig,psfrag,graphicx}
\usepackage{graphicx,psfrag}
\newtheorem{theorem}{Theorem}[section]
\newtheorem{lemma}[theorem]{Lemma}
\newtheorem{corollary}[theorem]{Corollary}
\newtheorem*{theoremA}{Main Theorem}
\newtheorem*{remark}{Remark}
\theoremstyle{definition}
\newtheorem{definition}[theorem]{Definition}
\numberwithin{equation}{section}
\newcommand{\N}{{\mathbb{N}}}
\newcommand{\Z}{{\mathbb{Z}}}

\newcommand{\dom}{\mathop{\mathrm{Dom}}}
\title[Transitive circle exchange transformations with flips]{Transitive circle exchange transformations\\ with flips}

\date{September 10, 2008}

\begin{document}

\title[Transitive circle exchange maps with flips]{Transitive circle exchange maps with flips}

\author{C.~Gutierrez}
\address{Instituto de Ci\^encias Matem\'aticas e de Computa\c c\~ao, Universidade de S\~ao Paulo, S\~ao Carlos - SP, BRAZIL}
\email{gutp@icmc.usp.br}

\author{S.~Lloyd}
\address{School of Mathematics and Statistics, University of New South Wales, Sydney NSW 2052, Australia}
\email{s.lloyd@unsw.edu.au}

\author{V.~Medvedev}
\address{Deptartment of Differential Equations, Institute of Applied
Mathematics and Cybernetics, Nizhny Novgorod State University, RUSSIA}
\email{medvedev@uic.nnov.ru}

\author{B.~Pires}
\address{Departamento de F\'isica e Matem\'atica, Faculdade de Filosofia, Ci\^encias e Letras de Ribeir\~ao Preto, Universidade de S\~ao Paulo, Ribeir\~ao Preto - SP, BRAZIL}
\email{benito@ffclrp.usp.br}

\author{E.~Zhuzhoma}
\address{Deptartment of Mathematics and Physics, Nizhny Novgorod State
Pedagogical University, RUSSIA}
\email{zhuzhoma@mail.ru}

\date{September 10, 2008}

\begin{abstract} We study the existence of transitive exchange maps with flips defined on the unit
circle $S^1$. We provide a complete answer to the question of whether there exists a transitive 
exchange map of $S^1$ defined on $n$ subintervals and having $f$ flips.
\end{abstract}
\maketitle
\section{Introduction}
Circle exchange transformations (CETs) and interval exchange transformations (IETs) form a 
fundamental class of objects studied in measurable dynamics, noted for their simple definition 
yet rich variety of behaviours (see \cite{Vi}). They arise naturally in the study of polygonal billiards (see \cite{BKM}), measured
foliations (see \cite{D-N,Ma}) and flat surfaces (see \cite{Zo1,Zo2}). 

Moreover, transitive CETs arise in the study of surface flows, as a factor of the Poincar\'e return map to a transversally embedded circle. 
The much-studied orientation-preserving CETs codify the
non-trivial recurrence of flows on orientable surfaces (see \cite{G1,Lev}). 
Advancement in the understanding of flows on \emph{non-orientable} surfaces, such as an answer for the open question
about the $C^r$-density, $r\ge 2$, of the Morse--Smale vector fields on non-orientable surfaces (see \cite{GB, Pe}),
requires study of the non-trivial recurrence in CETs that are not orientation-preserving, but rather reverse orientation on one or more subintervals, called ``flips''.

Generically, exchange maps with flips have periodic points and so are not transitive (see \cite{No2}). Nevertheless, for every even
$n\ge 4$ there exist transitive CETs with flips defined on $n$ subintervals (see \cite{G2,No1}).
In this article we extend this result: we determine for which pairs $(n,f)\in\mathbb{N}^2$ there exists a transitive CET defined on $n$ subintervals and having $f$ flips. 

We look for the least number of subintervals and flips required for a CET to be transitive. As remarked by Keane \cite{Ke1}, every CET with flips defined on two subintervals has a periodic point, and thus is not transitive. Thus we begin by looking at CETs defined on three subintervals.
The examples presented in this paper then work as building blocks of transitive, uniquely ergodic CETs having arbitrary numbers of subintervals and flips.

\section{Statement of the results}

Let $S^1 := [0,1]/(0\sim 1)$ be a circle endowed with the
orientation induced via the natural projection $[0,1]\to S^1$. 
We denote points of $S^1$ as if they were points of
$[0,1]$. Let $I_1,I_2,\ldots,I_n$ be a collection of pairwise
disjoint open subintervals of $S^1$ whose closure covers $S^1$.
A \emph{circle exchange transformation of $n$ subintervals},
shortly an $n$-CET, is an injective map $T:\bigcup_{i=1}^n I_i\subset
S^1\to S^1$ which acts as an isometry on each interval $I_i$,
$1\le i\le n$, and whose domain $\dom(T)=\bigcup_{i=1}^n I_i$ 
cannot be continuously extended to a bigger open subset of
$S^1$. A CET is called \emph{oriented} if it preserves the
orientation of all the domain subintervals, otherwise, if it reverses
the orientation of one or more domain subintervals, called \emph{flips},
it is said to be a CET \emph{with flips}. By replacing the circle
$S^1$ by a closed interval $[a,b]$, $0\leq a <b$, we get the definition 
of \emph{interval exchange transformation}. We shall use the expression 
\emph{exchange transformation} to refer to both CETs and IETs.
Notice that an $n$-CET (respectively, an $n$-IET) has necessarily $n$ (respectively, $n-1$)
discontinuity points. 

Let $\Gamma\in\{S^1,[a,b]\}$ and let $T:\Gamma\to\Gamma$ be an exchange map. The orbit of $p\in\Gamma$ is the set
$$O(p)=\{T^n(p)\mid n\in\mathbb{Z}\ {\rm and}\ p\in\dom(T^n)\}.$$
We say that $T$ is \emph{transitive} if there exists an orbit of $T$ that is dense in $\Gamma$. We say that the orbit of $p\in\Gamma$ is \emph{finite} if the set $O(p)$ has finite cardinality. Accordingly,
a point $p\in \Gamma-(\dom(T)\cup \dom(T^{-1}))$ has the finite orbit $O(p)=\{p\}$. 
A transitive exchange transformation is \emph{minimal} if it has no finite orbits.
The exchange transformation $T$ is called 
\emph{uniquely ergodic} if the Lebesgue measure on $\Gamma$ is the only
(up to scalar multiple) Borel measure invariant by $T$.
Notice that for exchange transformations, unique ergodicity implies transitivity.

The results of this paper may be combined to make the following statement.

\begin{theoremA}
Given $n\geq f\geq 1$, there exists a transitive CET of $n$ subintervals
with $f$ flips if and only if $n+f\geq 5$.
\end{theoremA}

We denote by $\mathcal{C}(n,f)$ the set of transitive CETs of $n$ subintervals with $f$ flips.
Since $C(2,f)$ is empty for $f=1,2$ (see \cite{Ke1}), it suffices to consider the case of $n\geq 3$. Thus the above result follows from the three statements:
\begin{itemize}
 \item [(i)] $\mathcal{C}(3,1)=\emptyset$ (Theorem \ref{31});
 \item [(ii)] $\mathcal{C}(3,f)\neq\emptyset$ for $f=2,3$ (Theorems \ref{32} and \ref{33});
 \item [(iii)] $\mathcal{C}(n,f)\neq\emptyset$ for every $n\ge 4$ and $1\le f\le n$ (Theorem \ref{nf}).
\end{itemize}

\section{$3$-CET with exactly one flip}\label{s:CETs}
The main result of this section is the
following theorem.

\begin{theorem}\label{31}
Every $3$-CET with exactly one flip has a periodic point and thus \mbox{$\mathcal{C}(3,1)=\emptyset$}.
\end{theorem}

The proof of the Theorem \ref{31} follows from the lemmas below. 

\begin{definition} Given $p\in S^1$, we identify $S^1-\{p\}$ with a
subinterval of $\mathbb{R}$ endowed with the total order induced by the
cyclic order. Accor\-dingly, given a subinterval $J$ of $S^1-\{p\}$
and $\delta\in\mathbb{R}$, with $\vert \delta\vert$ small enough, we may set
\mbox{$J+\delta=\{x+\delta:x\in J\}$} to be the shift of $J$
\mbox{by $\delta$}. Given subintervals $J_1,J_2$ of $S^1-\{p\}$, we
say that $J_1<J_2$ if $x<y$ for all $(x,y)\in J_1\times J_2$.
Analogously, if $\{J_i\}_{i=1}^r$, $r\ge 3$, is a collection of
subintervals of $S^1$, we say that \mbox{$J_1<J_2<\ldots<J_r$} if
every $(x_i)_{i=1}^r\in\Pi_{i=1}^r J_i$ is cyclically ordered.
\end{definition}

Throughout this section, we let $T:S^1\to S^1$ be a $3$-CET with one
flip defined on the domain $\dom(T)=S^1-\{a_0,a_1,a_2, a_3=a_0\}$, where
\mbox{$0=a_0<a_1<a_2<a_3=1$} and $I_i=({a_{i-1},a_i})$, $i=1,2,3$.
Without loss of generality, we may assume that $I_1$ is the flip of
$T$ and that \mbox{$T(I_1)<T(I_3)<T(I_2)$} (otherwise $T$ would be a
$2$-CET with flips). We denote by $T^0$ the identity map on $S^1$
and we suppose by contradiction that $T$ has no periodic points. We
say that $p\in S^1$ is a \emph{flipped} periodic point of $T$ if for some
$n\in\mathbb{N}$, $T^n(p)=p$ and $DT^n(p)=-1$.

\begin{lemma}\label{lf} There exists $N\geq 2$ such that ${T^n(I_1)}\cap {
I_1}=\emptyset$ for all $1\le n\le N-1$, and $T^N(I_1)\cap
I_1\neq\emptyset$. Furthermore, $\{T^n(I_1)\}_{n=0}^{N-1}$ are
pairwise disjoint.
\end{lemma}
\begin{proof}
The proof follows by induction. Firstly note that $T(I_1)\cap
I_1=\emptyset$, otherwise $T$ would have a flipped periodic point in
$I_1$. Now assume that for some $k\in\{1,2,\ldots\}$,
$\{T^n(I_1)\}_{n=0}^k$ are pairwise disjoint and that
$T^{k+1}(I_1)\cap I_1=\emptyset$. We claim that
$\{T^n(I_1)\}_{n=0}^{k+1}$ are pairwise disjoint. In fact, if
this was not the case then there would exist integers $1\le
m_1<m_2\le k+1$ such that \mbox{$T^{m_2}(I_1)\cap
T^{m_1}(I_1)\neq\emptyset$.} Hence, there would exist $x_2\in
T^{m_2-1}(I_1)$ and $x_1\in T^{m_1-1}(I_1)$ such that
$T(x_2)=T(x_1)$. By hypothesis, $T^{m_2-1}(I_1)\cap
T^{m_1-1}(I_1)=\emptyset$ and so $x_2\neq x_1$, which contradicts
the \mbox{injectivity} of $T$. This proves the claim. To finish the
proof, note that if $\{T^n(I_1)\}_{n=0}^k$ are pairwise disjoint
then $\mu(\bigcup_{n=0}^k\ T^n(I_1))=(k+1)\mu(I_1)$, where $\mu$
stands for the Lebesgue measure on $S^1$, which is invariant by $T$.
Hence, for some $k>0$, $T^k(I_1)\cap I_1\neq\emptyset$. We set $N$
to be the least $k$ with such property.
\end{proof}

By Lemma \ref{lf}, $D=\big(\bigcup_{n=0}^{N-1}
\,T^{n}(I_1)\big)\cap\{a_2\}$ has at most one point, i.e., $a_2$
belongs to at most one interval from $\bigcup_{n=0}^{N-1}
\,T^{n}(I_1)$. Hence, there are two cases to consider: either
$D=\emptyset$ or $D$ is a one-point set. Assume first that
$D=\emptyset$. Then $T^n({I_1})\cap \{a_0,a_1,a_2\}=\emptyset$ for
all $0\le n\le N-1$. Consequently, $I_1$ and $T^N(I_1)$ are
non--disjoint intervals with opposed orientations, and so $T$ has a
flipped periodic point in $I_1$, which contradicts the initial
assumption that $T$ has no periodic points. Therefore, we may assume
that $D\neq\emptyset$. Hence, for some $1\le n_0\le N-1$,
$T^{n_0}(I_1)\cap\{a_2\}\neq\emptyset$. Let $d\in I_1$ be such
$T^{n_0}(d)=a_2$ and let $I_1-\{d\}=I_\ell\cup I_r$, where
$I_\ell<I_r$ in $S^1-\{a_2\}$. Keeping these assumptions, we have
the following lemma.

\begin{lemma}\label{il} Either $a_0\in {T^{N}(I_\ell)}$ or $a_1\in {T^{N}(I_r)}$.
\end{lemma}
\begin{proof} Since $T^{N}(I_1)\cap I_1\neq\emptyset$, we have that either
$T^{N}(I_\ell)\cap I_1\neq\emptyset$ or $T^{N}(I_r)\cap
I_1\neq\emptyset$. We first assume that $T^{N}(I_r)\cap
I_1\neq\emptyset$. Note that each $J_n=T^n({I_r})$, $1\le n \le N$,
is an interval. For each $1\le n\le N$, let $k_n\ge 0$ and $1\le
q_n\le n_0$ be integers such that $n=k_n n_0+q_n$. We claim that for
all $1\le n\le N$,
$$J_n\subset S^1-\{a_0\} \quad\text{and}\quad
J_n=J_{q_n}-k_n\mu(I_r).$$ The claim holds trivially for all $1\le
n\le n_0$. Since $T(I_2)<T(I_1)<T(I_3)$ and
\mbox{$J_{n_0}=(c_{n_0},a_2)\subset I_2$,} we have that
$J_{n_0+1}=J_1-\mu(I_r)$. 

Now, given \mbox{$n_0+1< n\le N-1$,}
assume that the claim holds for all $1\le k\le n$. Let us prove that
it also holds for $n+1$. By the hypothesis of induction,
$J_{n-n_0}=J_{q_n}-{(k_n-1)}\mu(I_r)$ and \mbox{$J_n=J_{q_n}-k_n
\mu(I_r)$.} Hence, $J_n=J_{n-n_0}-\mu(I_r)$. Consequently (see
Figure \ref{figa}),
\begin{equation}\label{pe1}
 J_{n+1}=T(J_n)=T\big(J_{n-n_0}-\mu(I_r)\big)
=T(J_{n-n_0})-\mu(I_r)=J_{n-n_0+1}-\mu(I_r).
\end{equation}
 Using
the hypothesis of induction once more, we obtain that

\begin{equation}\label{pe2}
 J_{n-n_0+1}= \left\{ \begin{array}{ll}
         J_{q_n+1}-(k_n-1)\mu(I_r) & \mbox{if $1\le q_n<n_0$};\\
         J_{1}-k_n\mu(I_r)  & \mbox{if $q_n=n_0$}.\end{array} \right.
\end{equation}

Putting together the equations $(\ref{pe1})$ and $(\ref{pe2})$ we
reach
\begin{equation*}
 J_{n+1}= \left\{ \begin{array}{ll}
         J_{q_n+1}-k_n\mu(I_r)=J_{q_{n+1}}-k_{n+1}\mu(I_r) & \mbox{if $1\le q_n<n_0$};\\
         J_{1}-(k_n+1)=J_{q_{n+1}}-k_{n+1}\mu(I_r)  & \mbox{if $q_n=n_0$}.\end{array} \right.
\end{equation*}
which proves the claim.

\begin{figure}[ht]
  \centering
  \includegraphics[width=1.00\textwidth]{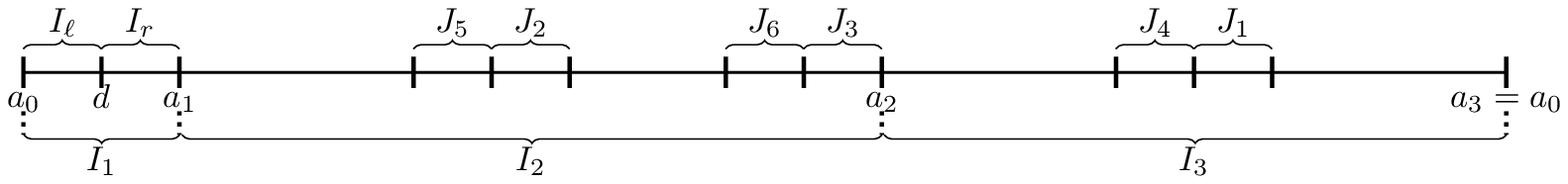}\\
  \caption{The iterates of $I_r$ under the map $T$.}\label{figa}
\end{figure}

Hence, we have that either $T^N(I_r)=I_r$ (which contradicts the
assumption that $T$ has no periodic points) or $a_1\in T^N(I_r)$.
The proof of the other case is similar.
\end{proof}

\begin{proof}[Proof of Theorem \ref{31}] It follows at once from Lemma
\ref{il}. In fact, the orientation of $T^{N}(I_\ell)$ and
$T^{N}(I_r)$ is opposed to the orientation of $I_\ell$ and $I_r$. By
Lemma \ref{il}, either $T^N(I_\ell)\cap I_\ell\neq\emptyset$ or
$T^N(I_r)\cap I_r\neq\emptyset$, which contradicts the initial
assumption that $T$ has no periodic points. This finishes the proof
because if a CET has a periodic point then it has a whole interval
of periodic points, and as such cannot be transitive. 
\end{proof}

\section{Rauzy induction for IETs with flips}

The Rauzy induction is a renormalisation algorithm for
interval exchange transformations \cite{Ra}. Using this method,
Veech \cite{Ve1} proved that almost all oriented IETs are uniquely ergodic. For explanations of the
oriented IET case, see \cite{Ve1}. Here we present the
Rauzy procedure for IETs with flips, see \cite{No2}. We
first introduce coordinates into the space of interval exchange
transformations.

Let $S_n$ be the permutation group on $\{1,2,\ldots,n\}$. We may represent any permutation $\pi\in S_n$ 
by the $n$-tuple
$$\pi=(\pi(1),\pi(2),\ldots,\pi(n))=(\pi_1,\pi_2,\ldots,\pi_n).
$$
We shall denote by $S_n^*$ the subset of $S_n$ made up by the
irreducible permutations: $\pi\in S_n^*$ if
$\pi(\{1,2,\ldots,k\})=\{1,2,\ldots,k\}$ implies $k=n$. Given $\theta\in \{-1,1\}^n$ and $\pi\in S_n$,
we let $\theta\,\pi=(\theta_1\pi_1 ,\theta_2\pi_2,\ldots,\theta_n\pi_n)$ be the signed
permutation obtained by componentwise multiplication of $\theta$ and $\pi$. We shall denote by
$\Sigma_n$ the set of signed permutations
$$
\Sigma_n=\{\theta\,\pi \mid \theta\in \{-1,1\}^n \ {\rm and}\ \pi\in S_n \}.
$$
We say that a signed permutation $p=\theta\,\pi$ is \emph{irreducible} if $\pi$ is irreducible. We shall denote by $\Sigma_n^*$ the set
of irreducible signed permutations on $n$ symbols.

To each $n$-IET $T:[0,\ell]\to [0,\ell]$, $\ell>0$, there corresponds a unique point
$(\lambda,p)$ in the space $\Lambda_n\times \Sigma_n$, where
$$\Lambda_n=\{\lambda=(\lambda_1,\ldots,\lambda_n)\mid\lambda_i>0,\,\forall i\}$$
is the positive cone. We say that
$(I_1,\ldots,I_n)$ is a \emph{partition} of $[0,\ell]$ if there exist
real numbers $0=a_0<a_1<\cdots<a_n=\ell$ such that
$I_i=(a_{i-1},a_i)$, $i=1,2,\ldots,n$. We let \mbox{$\vert J\vert=
\sup\,\{y-x\mid x,y\in J\}$.} Given an $n$-IET \mbox{$T:[0,\ell]\to [0,\ell]$},
we denote by $(I_1,\ldots,I_n)$ and $(J_1,\ldots,J_n)$,
respectively, the domain partition and the range partition defined
by $T$. We assign the data $(\lambda,p)\in \Lambda_n\times
\Sigma_n$ to $ T$ in the following way: for each
$i=1,2,\ldots,n$, $\lambda_i=\vert I_i\vert$, and $p=\theta\,\pi$, where 
 $\pi(i)$ is such that
$J_{\pi(i)}=T(I_{i})$, and $\theta_i=1$ (respectively $\theta_i=-1$) if $T$
preserves (respectively reverses) the orientation of $I_i$. Conversely, to
each $(\lambda,p)\in\Lambda_n\times\Sigma_n$,
there corresponds a unique IET of $k\le n$
subintervals that we shall denote by
$T=\mathcal{E}(\lambda,p)$.

Note that if $p\in
\Sigma_n^*$ is irreducible then any IET
$T=\mathcal{E}{(\lambda,p)}$ associated to the
signed permutation $p$ can be decomposed into several IETs with
smaller number of subintervals each. Since we are interested in
transitive exchange transformations, we shall only consider
irreducible permutations.

The \emph{Rauzy induction}
is the operator on the space of IETs that associates to
$T=\mathcal{E}(\lambda,p)$ the interval exchange
transformation $T'=\mathcal{E}(\lambda',p')$ which is the
Poincar\'e first return map induced by $T$ on the
subinterval $I(\lambda)=[0,\nu]$, where
$\vert\lambda\vert=\sum\lambda_i$, $p=\theta\,\pi$ and
\begin{eqnarray*}
\nu=\left\{
\begin{array}{ll}
\vert\lambda\vert-\lambda_{ \pi^{-1}(n)}, &\text{if}\ \lambda_{ \pi^{-1}(n)}<\lambda_n\\
\vert\lambda\vert-\lambda_n,              &\text{if}\ \lambda_{ \pi^{-1}(n)}>\lambda_n\\
\end{array}
\right..
\end{eqnarray*}
In this way, we obtain the map $\mathcal{T}_0:\Lambda_n\times
\Sigma_n\to \Lambda_n\times \Sigma_n$, where
$$\mathcal{T}_0(\lambda,p)=(\lambda',p')$$
and an
associated projective map
$$\mathcal{T}:\Delta_{n-1}\times \Sigma_n\to \Delta_{n-1}\times \Sigma_n,\quad
\Delta_{n-1}=\{\lambda\in\Lambda_n\mid\vert\lambda\vert=1\},$$ where
$$\mathcal{T}(\lambda,p)=\left(\frac{\lambda'}{\vert\lambda'\vert},p'\right).$$
Notice that the operators $\mathcal{T}_0$ and $\mathcal{T}$ are
defined only on the sets
\begin{eqnarray*}
\dom(\mathcal{T}_0)&=&\{(\lambda,\theta\,\pi)\in\Lambda_n\times
\Sigma_n\mid \lambda_n\neq\lambda_{ \pi^{-1}(n)}\}, \\
\dom(\mathcal{T})&=&\{(\lambda,\theta\,\pi)\in\Delta_{n-1}\times
\Sigma_n\mid \lambda_n\neq\lambda_{ \pi^{-1}(n)}\}.
\end{eqnarray*}

The operator $\mathcal{T}$ may be written in the form
\begin{eqnarray*}
\mathcal{T}(\lambda,\theta\,\pi)=\left\{
\begin{array}{lll}
\left(\frac{M_a(\theta\,\pi)^{-1}\lambda}{\vert M_a(\theta\,\pi)^{-1}\lambda\vert},a(\theta\,\pi)\right)&\text{if}&\  \lambda_{ \pi^{-1}(n)}<\lambda_n\\
\left(\frac{M_b(\theta\,\pi)^{-1}\lambda}{\vert M_b(\theta\,\pi)^{-1}\lambda\vert},b(\theta\,\pi)\right)&\text{if}&\  \lambda_{ \pi^{-1}(n)}>\lambda_n\\
\end{array}
\right.,
\end{eqnarray*}
where the transition matrices $M_a(\theta\,\pi),M_b(\theta\,\pi)\in
{\rm SL}_n(\mathbb{Z})$ and the transition maps $a,b:\Sigma_n\to
\Sigma_n$ are described below. For $i,j=1,\ldots,n$, denote
by $E_{ij}$ the $n\times n$ matrix of which the $(i,j)$th element is
equal to $1$, and all the others elements are equal to $0$. Let $E$
be the $n\times n$ identity matrix. The transition matrices are
defined by:
\begin{eqnarray*}
M_a(\theta\,\pi)&=&E+E_{n, \pi^{-1}(n)},\\
M_b(\theta\,\pi)&=&\sum_{i=1}^{ \pi^{-1}(n)}
E_{i,i}+E_{n,s(\theta\,\pi)}+\sum_{i= \pi^{-1}(n)}^{n-1} E_{i,i+1},
\end{eqnarray*}
where $s(\theta\,\pi)= \pi^{-1}(n)+(1+\theta_{ \pi^{-1}(n)})/2$.

We now define the transition maps. When $\theta_n=+1$ the transition map
$a$ is defined by
$$a(\theta\,\pi)_i=\left\{
\begin{array}{ll}
\theta_i \pi_i, & \pi_i\leq \pi_n, \\
\theta_i (\pi_n+1), & \pi_i=n, \\
\theta_i (\pi_i+1), & \textrm{otherwise},\\
\end{array}
\right.
$$
and when $\theta_n=-1$, we have
$$
a(\theta\,\pi)_i=\left\{
\begin{array}{ll}
\theta_i \pi_i, & \pi_i\leq \pi_n-1, \\
-\theta_i \pi_n, & \pi_i=n, \\
\theta_i (\pi_i+1), & \textrm{otherwise}.\ \end{array} \right.
$$
The transition map $b$ when $\theta_{\pi^{-1}(n)}=+1$ is defined by
$$
b(\theta\,\pi)_i=\left\{
\begin{array}{ll}
\theta_i \pi_i, & i\leq  \pi^{-1}(n), \\
\theta_n \pi_n, & i=  \pi^{-1}(n)+1, \\
\theta_{i-1} \pi_{i-1}, & \textrm{otherwise},\\
\end{array}
\right.
$$
and when $\theta_{\pi^{-1}(n)}=-1$ by
$$
b(\theta\,\pi)_i=\left\{
\begin{array}{ll}
\theta_i \pi_i, & i\leq  \pi^{-1}(n)-1, \\
-\theta_n\pi_n, & i=  \pi^{-1}(n), \\
\theta_{i-1} \pi_{i-1}, & \textrm{otherwise}.\\
\end{array}
\right.
$$

\section{Renormalisable interval exchange maps}

We shall denote by $\mathcal{G}_n$ the graph with vertices in $\Sigma_n^*$ and edges in $\mathbb{E}_n=\mathbb{E}_n^a\cup \mathbb{E}_n^b$, where
\begin{eqnarray*}
{\mathbb{E}}_n^a&=&\{(p,q)\in\Sigma_n^*\times\Sigma_n^*\mid q=a(p)\}\\
{\mathbb{E}}_n^b&=&\{(p,q)\in\Sigma_n^*\times\Sigma_n^*\mid q=b(p)\}
\end{eqnarray*}
and $a,b:\Sigma_n\to\Sigma_n$ are the transitions maps defined previously. Since ${\mathbb{E}}_n^a\cap {\mathbb{E}}_n^b=\emptyset$, there is a function $\tau:\mathbb{E}_n\to\{a,b\}$ which assigns to each edge $e$ of $\mathcal{G}_n$ its type: \mbox{$a$ (if $e\in {\mathbb{E}}_n^a$)} or $b$ (if
$e\in {\mathbb{E}}_n^b$). Let $N$ be a positive integer. A (directed) \emph{finite path} of $\mathcal{G}_n$
is a sequence $\{p^{(i)}\}_{i=0}^N$ such that $(p^{(i-1)},p^{(i)})\in \mathbb{E}_n$ for all $1\le i\le N$. To each finite path
$\{p^{(i)}\}_{i=0}^N$, there are associated a \emph{type itinerary} $\{t_i\}_{i=1}^N$, defined by $t_i=\tau(p^{(i-1)},p^{(i)})$
for $i=1,2,\ldots,N$, and a \emph{matrix itinerary} $\{M_{\gamma}^{(i)}\}_{i=1}^{N}$ consisting of the product
of transition matrices
\begin{equation*}
M_\gamma^{(i)}=M_{t_1}(p^{(0)})\cdot M_{t_2}(p^{(1)})\cdots M_{t_i}(p^{(i-1)}). 
\end{equation*}
By replacing the set $\{0,1,\ldots,N\}$ by $\N\cup\{0\}$, we obtain the definition of \emph{infinite path}. In the same way as above, we get
the definition of type itinerary $\{t_i\}_{i=1}^\infty$ and of matrix itinerary $\{M_{\gamma}^{(i)}\}_{i=1}^{\infty}$
associated to an infinite path $\{p^{(i)}\}_{i=0}^\infty$. We say that a finite path $\{p^{(i)}\}_{i=0}^{N}$ is a \emph{cycle} if $p^{(N)}=p^{(0)}$.  To every
cycle $\gamma=\{p^{(i)}\}_{i=0}^{N}$ there are associated the \emph{product matrix} $M_\gamma^{(N)}$ and the infinite path
$\{p^{(i)}\}_{i=0}^\infty$ where $p^{(j)}:=p^{(j\,{\rm mod}\,N)}$ if $j>N$. We call such infinite path a \emph{periodic path
of period $N$}. 

We say that a matrix
$A\in{\rm SL}_n(\Z)$ is \emph{eventually positive} if $A^n$ is a positive matrix for some positive integer $n$. A
 cycle $\gamma=\{p^{(i)}\}_{i=0}^{N}$ is \emph{special} if the associated product matrix $M=M^{(N)}_\gamma$ is
eventually positive. A path $\{p^{(i)}\}_{i=0}^N$ is a \emph{special path} if $p^{(N)}$ belongs to a special cycle. A periodic path $\gamma=\{p^{(i)}\}_{i=0}^\infty$  of period $N$ is \emph{special} if $M_\gamma^{(N)}$ is eventually positive.

We say that two cycles $\{p^{(i)}\}_{i=0}^{N}$ and $\{q^{(j)}\}_{j=0}^{N}$
are \emph{shift equivalent} if there exists $r\in\Z$ such that $q^{(i)}=p^{((i+r){\,{\rm mod}\, N})}$ for all
$i\in\{0,1,\ldots,N\}$. The property of being a special cycle is invariant under shift equivalence. 

We say that $T=\mathcal{E}(\lambda^{(0)},p^{(0)})$ is \emph{infinitely renormalisable} if
there exists a sequence of pairs $\{(\lambda^{(i)},p^{(i)})\}_{i=0}^\infty$ in 
$\dom(\mathcal{T}_0)\cap (\Lambda_n\times\Sigma_n^*)$ such that for all $i>0$,
\begin{equation*}\label{ir}
(\lambda^{(i)},p^{(i)})=
\mathcal{T}_0(\lambda^{(i-1)},p^{(i-1)}).
\end{equation*}
In this case, the infinite path $\{p^{(i)}\}_{i=0}^\infty$ is called the \emph{path
corresponding to $T$}. Not every infinite path in $\mathcal{G}_n$ corresponds
to an infinitely renormalisable IET. Nevertheless, special periodic paths (and so special cycles) always correspond to infinitely renormalisable
IETs. 

Given an infinite path $\gamma=\{p^{(i)}\}_{i=0}^\infty$ of $\mathcal{G}_n$, let
$$C(\gamma)=\bigcap_{i=1}^\infty M_\gamma^{(i)}{\Lambda_n}.$$
A vector $\lambda$ of $\Lambda_n$ is in $C(\gamma)$ if and only if $\mathcal{E}(\lambda,p^{(0)})$ is infinitely renormalisable and its corresponding path
is $\gamma$. As in the oriented case (see \cite{Vi}), if $T=\mathcal{E}(\lambda,p^{(0)})$ is infinitely renormalisable,
 then it is minimal. The dimension of $\mathcal{C}(\gamma)$ is strictly less than $n$, and $C(\gamma)$ is linearly
isomorphic to $\mathcal{M}(T)$, the set of invariant measures of $T$
that are positive on open sets. Thus $T$ is uniquely ergodic if and
only if $\mathcal{C}(\gamma)$ is a half-line.

We now explain how to construct infinitely renormalisable IETs with flips. Let $\{p^{(i)}\}_{i=0}^N$ be a special cycle and let
$\gamma=\{p^{(i)}\}_{i=0}^\infty$ be the corresponding special periodic path. By hypothesis the matrix
$M=M_\gamma^{(N)}$ is eventually positive.
The Perron--Frobenius Theorem asserts that $M$ has a
positive, simple eigenvalue $\sigma$ of maximum  modulus, and a
unique eigenvector $\lambda\in\Delta_{n-1}$ such that
$M\lambda=\sigma\lambda$. Thus
$$
\lambda\in\mathcal{C}(\gamma)=\bigcap_{i=1}^\infty
M^{(i)}_\gamma\Lambda_n\subset \bigcap_{j=1}^\infty
M^j\Lambda_n=\textrm{span}(\lambda).
$$
Hence, taking $\lambda^{(0)}=\lambda$, we have that the domain of
$\mathcal{E}(\mathcal{T}_0^{N}(\lambda^{(0)},p^{(0)}))$
is $[0,1/\sigma]$, and
$(\lambda^{(0)},p^{(0)})$ is a
periodic point of $\mathcal{T}$ of period $N$. Thus
$T_0=\mathcal{E}(\lambda^{(0)},p^{(0)})$ is
infinitely renormalisable, and therefore minimal. By the uniqueness
of the positive eigendirection of $M$, we have that
$\mathcal{C}(\gamma_0)$ is a half-line, and so $T_0$ is uniquely
ergodic. As $(\lambda^{(0)},p^{(0)})$ is a fixed
point of the renormalisation operator $\mathcal{T}^N$, we say that
the IET
$T_0=\mathcal{E}(\lambda^{(0)},p^{(0)})$ is {\it
self-induced}.

\section{Minimal $3$-CET with two flips}\label{s:3-2}

In this section we present a special cycle of $\mathcal{G}_4$, from which we construct a
self-induced interval exchange transformation of four subintervals
with two flips. We also make use of this cycle in each of the
subsequent sections.

Given a pair $(n,f)\in\N^2$, we let $\mathcal{I}(n,f)$ (respectively, $\mathcal{C}(n,f)$) denote the set of transitive IETs
(respectively, CETs) of $n$ subintervals with $f$ flips.

\begin{theorem}\label{32}
There exists a minimal, uniquely ergodic
$4$-CET with two flips. In particular, $\mathcal{I}(4,2)\neq\emptyset$ and $\mathcal{C}(3,2)\neq\emptyset$.
\end{theorem}
\begin{proof} Let $\gamma_1=\{p^{(i)}\}_{i=0}^{9}$ be the
cycle of $\mathcal{G}_4$ that starts at the signed permutation
$p^{(0)}=(-3,+4,+1,-2)$, has type itinerary $\{a,b,a,b,b,a,b,a,b\}$ and is described in the Figure
\ref{cycle}.

\begin{figure}[ht]
  \centering
  \includegraphics[width=0.7\textwidth]{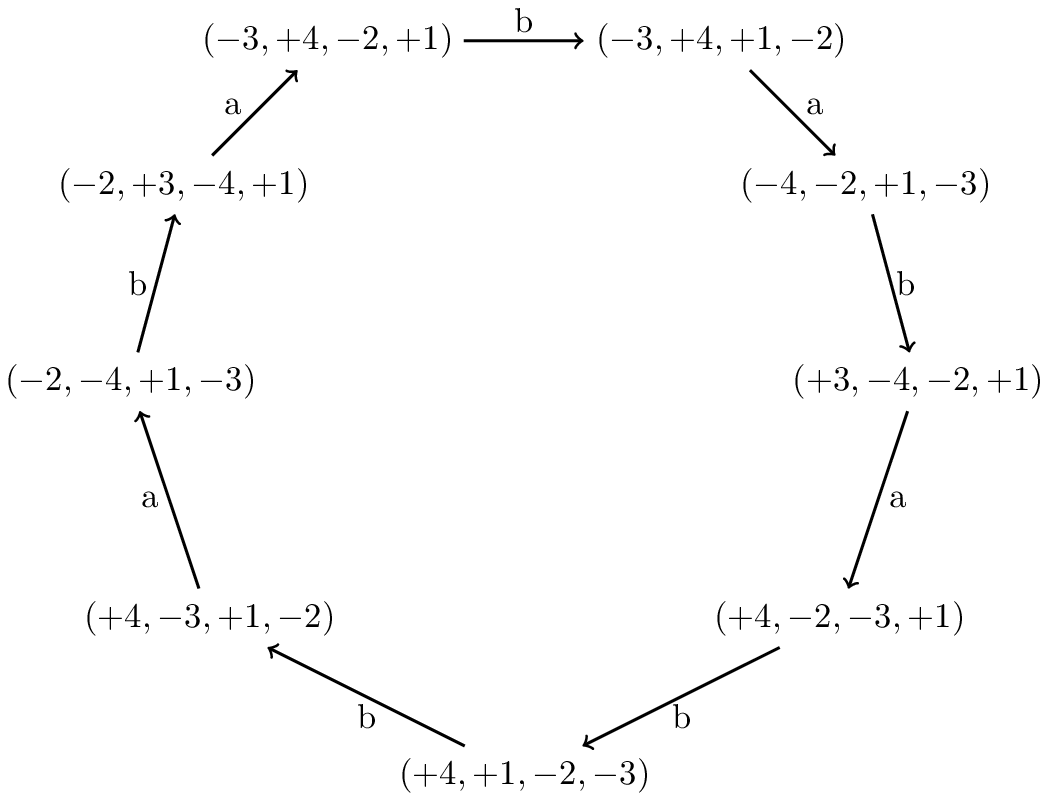}\\
  \caption{A special cycle of $\mathcal{G}_4$.}\label{cycle}
\end{figure}

The product matrix associated to the cycle $\gamma_1$ is

$$
M=M_{\gamma_1}^{(9)}=\left(
\begin{array}{cccc}
2 & 2 & 3 & 2 \\
0 & 1 & 1 & 1 \\
1 & 1 & 2 & 1 \\
1 & 2 & 3 & 3 \\
\end{array}
\right).
$$
The matrix $M$ is eventually positive and so $\gamma_1$ is a special cycle. If $\sigma$ is the Perron--Frobenius eigenvalue
of $M$ and $\lambda\in\Delta_{3}$ is the corresponding eigenvector then
the $4$-IET $T_1=\mathcal{E}(\lambda^{(0)},p^{(0)})$ of
$[0,1]$ with $\lambda^{(0)}=\lambda$ and permutation $p^{(0)}=(-3,+4,+1,-2)$
is self-induced on the interval $[0,1/\sigma]$ (see Figure
\ref{2flips}). Notice that the path of $\mathcal{G}_4$ corresponding to $T_1$ is the infinite periodic path
generated by $\gamma_1$, i.e., the path $\{p^{(i)}\}_{i=0}^\infty$, where $p^{(j)}:=p^{(j{\,\rm mod\,}(9))}$ for
all $j>9$.

If the endpoints of $[0,1]$ are identified, $T_1$ becomes a minimal,
uniquely ergodic $3$-CET with two flips as required because the
union $I_2\cup \{a_2\}\cup I_3$ is mapped onto $J_4\cup \{0\sim
1\}\cup J_1$. 
\end{proof}

\section{Minimal $3$-CET with three flips}\label{s:3-3} 

In this section we prove the following result:

\begin{theorem}\label{33}
There exists a minimal, uniquely ergodic
$4$-CET with three flips. In particular, $\mathcal{I}(4,4)\neq\emptyset$ and $\mathcal{C}(3,3)\neq\emptyset$.
\end{theorem}
\begin{proof} Let $\gamma_2=\{q^{(j)}\}_{j=0}^4$ be the path of $\mathcal{G}_4$ that starts at the signed permutation
$q^{(0)}=(-3,-1,-4,-2)$, has type itinerary $\{a,b,b,a\}$, and is described in \mbox{Figure \ref{path1}.} Note that $\gamma_2$ is a special path, since 
$q^{(4)}=(+4,-3,+1,-2)$ is a vertex in the special cycle $\gamma_1$.

\begin{figure}[ht]
    \centering
        \includegraphics[width=0.9\textwidth]{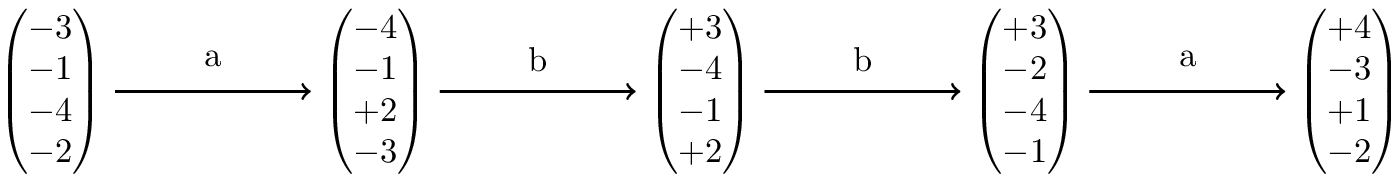}
    \caption{A path in $\mathcal{G}_4$ terminating on the cycle of Figure \ref{cycle}.}
    \label{path1}
\end{figure}

The product matrix for the path $\gamma_2$ is
$$
M^{(4)}_{\gamma_2}=\left(
\begin{array}{cccc}
1 & 1 & 1 & 0 \\
0 & 0 & 1 & 1 \\
0 & 1 & 0 & 0 \\
1 & 1 & 0 & 0 \\
\end{array}
\right).
$$
Setting $\mu:=M_{\gamma_2}^{(4)}(M_{\gamma_1}^{(5)})^{-1}\lambda$ yields a $4$-IET 
$T_2:=\mathcal{E}(\mu/\vert\mu\vert,q^{(0)})$ of $[0,1]$ with four flips (see Figure \ref{fig:3flips}) that is minimal and 
uniquely ergodic, since there is a subinterval to which the
Poincar\'e first return map is minimal and uniquely ergodic. Identifying the
endpoints of $[0,1]$, since the interval $I_2\cup \{a_2\}\cup I_3$
is mapped onto the interval $J_4\cup \{0\sim 1\}\cup J_1$, we have a minimal,
uniquely ergodic $3$-CET with three flips.
\end{proof}

\section{Minimal $4$-CET with one flip}

In this section we prove the following result, thus showing
that \mbox{Theorem \ref{31}} cannot be extended to $n$-CETs, $n>3$.

\begin{figure}[ht]
    \centering
    \psfrag{T1}{$T_1$}
    \psfrag{T2}{$T_2$}
    \psfrag{T3}{$T_3$}
        \includegraphics[width=\linewidth]{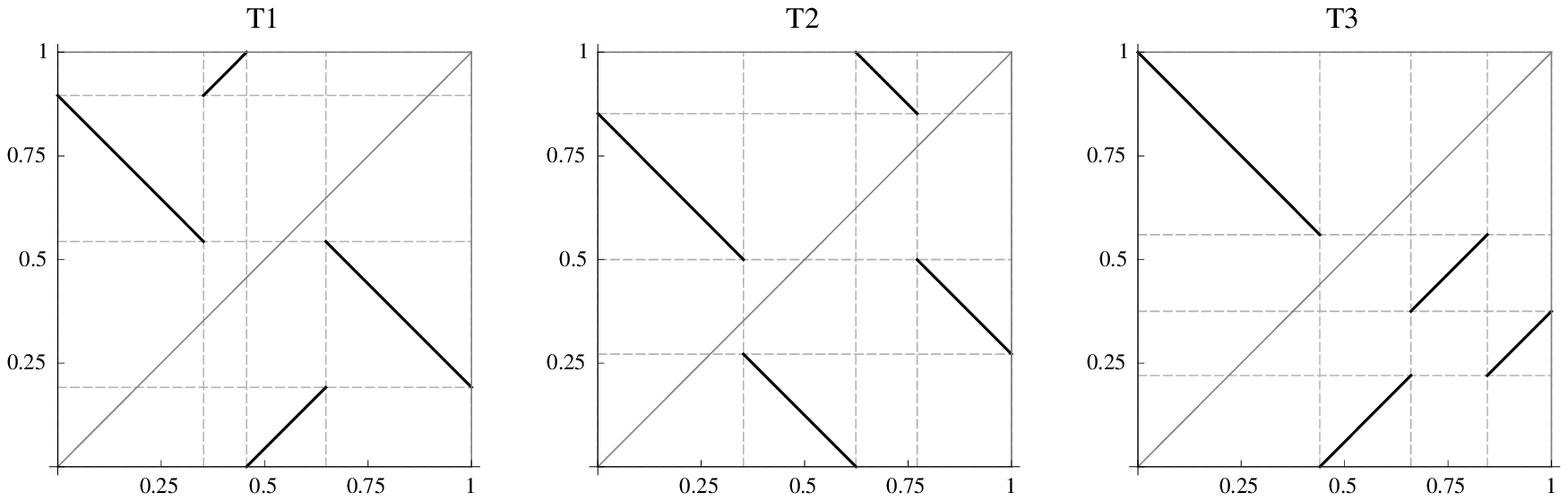}
    \caption{Graphs of the minimal IETs $T_1,T_2,T_3$.}
    \label{fig:1flip}
    \label{2flips}
    \label{fig:3flips}
\end{figure}

\begin{theorem}\label{41}
There exists a minimal, uniquely ergodic
$4$-IET with one flip. In particular, $\mathcal{I}(4,1)\neq\emptyset$ and $\mathcal{C}(4,1)\neq\emptyset$.
\end{theorem}
\begin{proof} We define $\gamma_3=\{r^{(k)}\}_{k=0}^4$ to be the path
in $\mathcal{G}_4$ that starts at $r^{(0)}=(-4,+1,+3,+2)$, has type itinerary $\{b,b,a,b\}$ and finishes
at the vertex $r^{(4)}=(-3,-1,-4,-2)$ of $\gamma_2$ (see Figure \ref{path2}).

\begin{figure}[ht]
    \centering
        \includegraphics[width=0.9\textwidth]{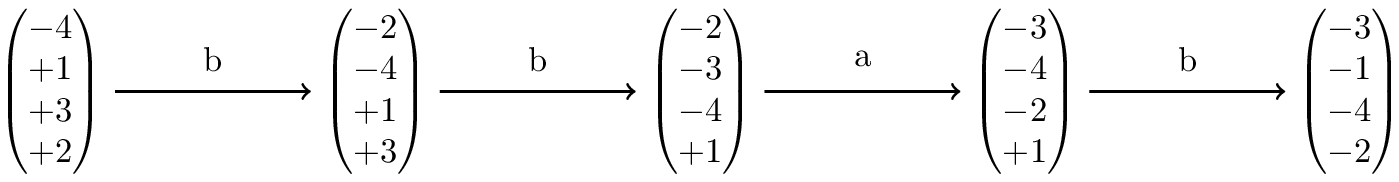}
   \caption{A path in $\mathcal{G}_4$ terminating at the start of the path in Figure \ref{path1}.}
    \label{path2}
\end{figure}

The product matrix for the path $\gamma_3$ is
$$
M^{(4)}_{\gamma_3}=\left(
\begin{array}{cccc}
1 & 1 & 1 & 1 \\
0 & 1 & 0 & 1 \\
0 & 1 & 1 & 0 \\
1 & 0 & 0 & 0 \\
\end{array}
\right).
$$
Multplying the vector $\mu/|\mu|$ by this matrix and
normalizing gives the required length vector $\xi$. Thus we have constructed a $4$-IET
$T_3:=\mathcal{E}(\xi,r ^{(0)})$ (see
Figure \ref{fig:1flip}) that has exactly one
flip, and since it has a subinterval to which the Poincar\'e first
return map is minimal and uniquely ergodic, is itself minimal and uniquely ergodic also.
\end{proof}

\section{Transitive CETs with arbitrary number of subintervals and flips}\label{s:arbf}

In this section we provide a proof of the following result.

\begin{theorem}\label{nf}
$\mathcal{C}(n,f)\neq\emptyset$ for all $n\ge 4$ and $1\le f\le n$.
\end{theorem}

\begin{definition}[fake discontinuity] We say
that an $n$-IET $T=\mathcal{E}(\lambda,\theta\,\pi)$ has a fake
discontinuity if precisely one of the following cases happens:
\begin{itemize}
\item[$(a)$] for exactly one $1\le i\le n-1$, $\pi(i)=n$,
$\pi(i+1)=1$, and $\theta_i=\theta_{i+1}=+1$,  
\item [$(b)$] for
exactly one $1\le i\le n-1$, $\pi(i)=1$, $\pi(i+1)=n$, and
$\theta_i=\theta_{i+1}=-1$,
\item[$(c)$] $\pi(1)=\pi (n)+1$ and $\theta_1=\theta_n=+1$,
\item[$(d)$] $\pi(1)=\pi(n)-1$ and $\theta_1=\theta_n=-1$.
\end{itemize}
\end{definition}

Henceforth, we call a transitive $n$-IET (respectively $n$-CET) with $f$ flips an $(n,f)$-IET (respectively
an $(n,f)$-CET). An $(n,f)$-IET without fake discontinuities is isomorphic to an $(n,f)$-CET whereas
an $(n,f)$-IET with one fake discontinuity is isomorphic to either an $(n-1,f)$-CET or to an $(n-1,f-1)$-CET.
Notice that an $(n,n)$-IET with one fake discontinuity is always isomorphic to an $(n-1,n-1)$-CET.

We shall need the following result.

\begin{lemma}\label{lem45} There exist the following kind of transitive IETs:
\begin{itemize}
 \item [(a)] A $(4,f)$-IET without fake discontinuities, for each $f=1,2,3$;
 \item [(b)] A $(5,5)$-IET with exactly one fake discontinuity;
 \item [(c)] A $(5,5)$-IET without fake discontinuities.
\end{itemize}
\end{lemma}
\begin{proof} (a) Theorem \ref{41} gives a (4,1)-IET without fake discontinuities. The permutation 
$p^{(2)}=(+3,-4,-2,+1)$ of $\gamma_1$ and $q^{(1)}=(-4,-1,+2,-3)$ of $\gamma_2$ provide, respectively, a $(4,2)$-IET and a $(4,3)$-IET, both
without fake discontinuities. As for items (b) and (c), let $\gamma_b=\{p_b^{(i)}\}_{i=0}^{12}$ be the path starting at $(-3,-1,-5,-2,-4)$ having type itinerary $\{b,a,a,b,b,b,a,b,a,b,a,b,a\}$ and
let $\gamma_c=\{p_c^{(j)}\}_{j=0}^{21}$ be the path starting at \mbox{$(-2,-3,-4,-5,-1)$} having type itinerary 
$$\{a,a,b,b,a,b,b,a,a,a,b,b,a,b,b,a,b,a,b,a,a,b\}.$$
We claim that these paths are special. Indeed, the fifth permutation of $\gamma_b$, $p_b^{(4)}=(4,-5,-1,-3,2)$, belongs to a cycle whose product matrix
$B$ is eventually positive and the tenth permutation of $\gamma_c$, $p_c^{(9)}=(-4,5,1,-2,-3)$, belongs to a cycle whose product matrix $C$ is eventually positive. We have

\begin{equation*}
\begin{array}{ccc}
B=\begin{pmatrix}
 1 & 1 & 0 & 0 & 0\\
 1 & 0 & 3 & 4 & 3\\
 1 & 1 & 1 & 1 & 1\\
 0 & 0 & 1 & 2 & 2\\
 0 & 0 & 1 & 1 & 0
\end{pmatrix},
&   & C=
\begin{pmatrix}
 2 & 2 & 3 & 4 & 3\\
 0 & 1 & 1 & 1 & 1\\
 1 & 1 & 2 & 1 & 1\\
 0 & 0 & 0 & 2 & 1\\
 1 & 2 & 3 & 3 & 3
\end{pmatrix}.
\end{array}
\end{equation*}
To obtain the required IETs we proceed as in the previous sections. Let $\lambda_b\in\Delta_{4}$ and $\lambda_c\in\Delta_{4}$ be the
Perron--Frobenius eigenvectors of $B$ and $C$, respectively. Set $\mu_b=M_{\gamma_b}^{(4)}\lambda_b$ and $\mu_c=M_{\gamma_c}^{(9)}\lambda_c$.
The minimal $5$-IET $T_b=\mathcal{E}(\mu_b/\mu_b,p_b^{(0)})$ has $5$ flips and one fake discontinuity whereas the minimal
$5-$IET $T_c=\mathcal{E}(\mu_c/\mu_c,p_c^{(0)})$ has $5$ flips and no fake discontinuities.
\end{proof}

\begin{corollary}\label{cor4} $\mathcal{C}(4,f)\neq\emptyset$, for all $f=1,2,3,4$.
\end{corollary}
\begin{proof}
 It follows from itens (a) and (b) of Lemma \ref{lem45}.
\end{proof}

Now we introduce two operators in the space of the interval exchange transformations that will allow us
to construct IETs with arbitrary numbers of subintervals and flips. Given an $(n,f)$-IET $T:[0,1]\to
[0,1]$, we can define an \mbox{$(n+1,f)$-IET}
$T_1:[0,2]\to [0,2]$ and an $(n+2,f+2)$-IET
$T_2:[0,3]\to [0,3]$ such that the Poincar\'e first return map
induced by $T_1$ (respectively $T_2$) on $[0,1]$ is $T$. We set:
\begin{equation*}
\begin{array}{cc}
T_1(x)=\left\{
 \begin{array}{ll}
T(x)+1 & \text{if}\ x\in [0,1]\\
x-1    & \text{if}\ x\in [1,2]
 \end{array}
\right., &
T_2(x)=\left\{
 \begin{array}{ll}
T(x)+1 &\text{if}\ x\in [0,1]\\
-x+4   &\text{if}\ x\in [1,2]\\
-x+3   &\text{if}\ x\in [2,3]\\
 \end{array}
\right..
\end{array}
\end{equation*}
Let $\rho:\Lambda_n\times \Sigma_n\to \Delta_{n-1}\times\Sigma_n$ be the projection map defined by
$\rho(\lambda,\theta\,\pi)=(\lambda/\vert\lambda\vert,\theta\,\pi)$
and let $\alpha$ and $\beta$ be the operators in the space of the
transitive interval exchange transformations which at a
$(n,f)$-IET $T:[0,1]\to [0,1]$ takes the value $\alpha(T)=\rho(T_1)$
and $\beta(T)=\rho(T_2)$, respectively. The operators
$\alpha$ and $\beta$ always produce IETs without fake
discontinuities which so correspond to CETs with the same number of
subintervals and flips. 

\begin{figure}[htbp]
    \centering
    		\psfrag{a}{$\alpha$}
    		\psfrag{b}{$\beta$}
    		\psfrag{1}{$1$}
				\psfrag{2}{$2$}
    		\psfrag{3}{$3$}
    		\psfrag{4}{$4$}
    		\psfrag{5}{$5$}
    		\psfrag{6}{$6$}
    		\psfrag{7}{$7$}
    		\psfrag{8}{$8$}
    		\psfrag{9}{$9$}
        \includegraphics[width=0.90\textwidth]{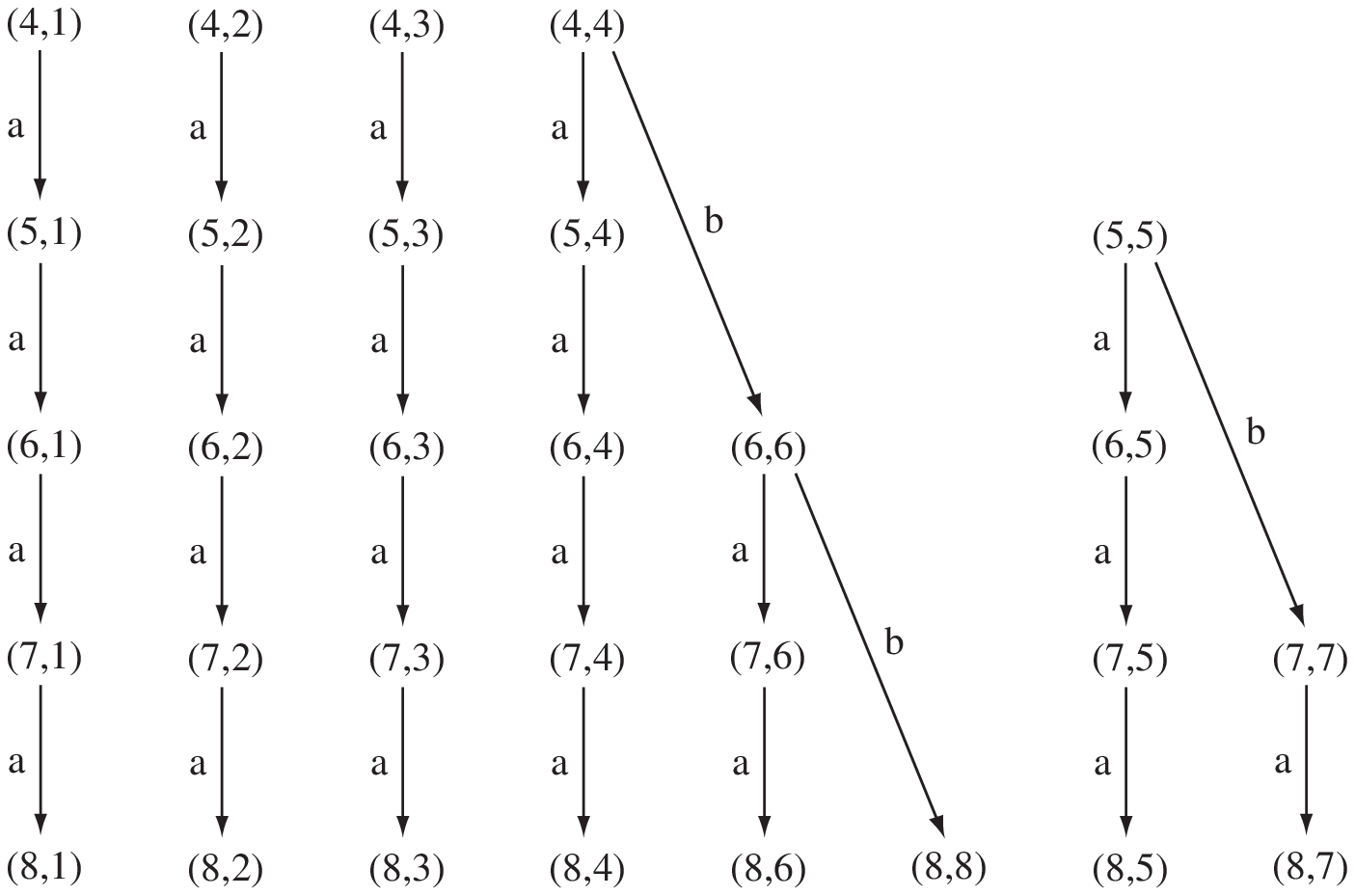}
    \caption{Tree of the $(n,f)$-IETs without fake discontinuities \mbox{generated} by $\alpha$
and $\beta$.}
\label{arvore}
\end{figure}

We indicate now how to proceed in order to extend the results of
this paper for more subintervals. The case of four subintervals was
settled by the Corollary \ref{cor4}. The case of five subintervals can
be approached in the following way. By Lemma \ref{lem45}, there
exists a $(4,i)$-IET $E_i:[0,1]\to [0,1]$, $i=1,2,3$. By  Theorem \ref{32}, we have a $(4,4)$-IET $E_4:[0,1]\to [0,1]$ (with a
fake discontinuity). The $(5,i)$-IETs $\alpha(E_i)$, $i=1,2,3,4,$
have no fake discontinuities and so are isomorphic to $(5,i)$-CETs.
The remaining case $(5,5)$ is covered by Lemma \ref{cor4}. The case
of $n>5$ subintervals can be obtained by repeated applications of the
operators $\alpha$ and $\beta$ to our $4$-IETs and $5$-IETs
examples, as it can be seen in the Figure \ref{arvore} (a pair
$(n,f)$ means a uniquely ergodic transitive $(n,f)$-IET).
Note that in Figure \ref{arvore}, the vertex $(4,4)$ has a fake 
discontinuity, but all the other ones are free of fake discontinuities. 
We have proved Theorem \ref{nf}.

\begin{remark} All the examples presented in this paper are uniquely ergodic. By applying the methods of Section \ref{s:arbf} to the non-uniquely ergodic oriented IET example by Keane \cite{Ke2}, one can easily construct 
transitive non-uniquely ergodic IETs with flips. However we leave open the problem of finding minimal non-uniquely ergodic IETs and CETs 
with flips.
\end{remark}

\medskip\noindent{\it Acknowledgments}. We would like to thank D.~Anosov and A.~Nogueira for their
helpful comments to a previous version of this article. This work began while the fifth author was
visiting University of S\~ao Paulo at S\~ao Carlos (ICMC-USP) in
March-April 2007. He would like to thank Carlos
Gutierrez for the hospitality, and the support of University of S\~ao Paulo
which made the visit possible.

VM and EZ were partially supported by RFFI grants 08-01-00547a and 05-01-00501. 
CG was partially supported by FAPESP Grant 03/03107-9 and CNPq Grants 470957/2006-9 and
306328/2006-2, Brazil. BP and SL were fully supported by FAPESP Grants 06/52650-5 and 06/60600-8, Brazil,
respectively.


\bibliographystyle{amsplain}

\end{document}